\documentclass[12pt,thmsa]{article}
\usepackage{amsfonts}

\input{tcilatex}
\begin{document}

\author{Karl-Georg Schlesinger \qquad \\
Erwin Schr\"{o}dinger Institute for Mathematical Physics\\
Boltzmanngasse 9\\
A-1090 Vienna, Austria\\
e-mail: kgschles@esi.ac.at}
\title{Some remarks on $q$-deformed multiple polylogarithms}
\date{}
\maketitle

\begin{abstract}
We introduce general $q$-deformed multiple polylogarithms which even in the
dilogarithm case differ slightly from the deformation usually discussed in
the literature. The merit of the deformation as suggested, here, is that $q$%
-deformed multiple polylogarithms define an algebra, then (as in the
undeformed case). For the special case of $q$-deformed multiple $\zeta $%
-values, we show that there exists even a noncommutative and
noncocommutative Hopf algebra structure which is a deformation of the
commutative Hopf algebra structure which one has in the classical case.
Finally, we discuss the possible correspondence between $q$-deformed
multiple polylogarithms and a noncommutative and noncocommutative self-dual
Hopf algebra recently introduced by the author as a quantum analog of the
Grothendieck-Teichm\"{u}ller group.
\end{abstract}

\section{Motivation}

In \cite{CL} a family of noncommutative deformations $S^{4,\theta }$ of the
four dimensional sphere has been introduced from the instanton algebra.
Concretely, the algebra of functions on $S^{4,\theta }$ is given by the
generators $t,\alpha ,\alpha ^{*},\beta ,\beta ^{*}$ and relations 
\begin{eqnarray*}
\alpha \alpha ^{*} &=&\alpha ^{*}\alpha ,\ \beta \beta ^{*}=\beta ^{*}\beta
\\
\alpha \beta &=&\lambda \beta \alpha ,\ \alpha ^{*}\beta =\overline{\lambda }%
\beta \alpha ^{*} \\
0 &=&\alpha \alpha ^{*}+\beta \beta ^{*}+t^2-t
\end{eqnarray*}
and the requirement that $t$ commutes with all the other generators. With $%
\varphi ,\psi $ angles, 
\[
0\leq \varphi \leq \frac \pi 2,\ -\frac \pi 2\leq \psi \leq \frac \pi 2 
\]
and $u,v$ the unitary generators of the algebra of smooth functions on the
noncommutative 2-torus, one has 
\begin{eqnarray*}
\alpha &=&\frac u2\cos \varphi \cos \psi \\
\beta &=&\frac v2\sin \varphi \cos \psi \\
t &=&\frac 12+\frac 12\sin \psi
\end{eqnarray*}
and the Dirac operator corresponding to the round metric on the sphere is
given as (see \cite{CL}) 
\begin{eqnarray*}
D &=&\left( \cos \varphi \cos \psi \right) ^{-1}\delta _1\gamma _1+\left(
\sin \varphi \cos \psi \right) ^{-1}\delta _2\gamma _2 \\
&&+\frac 1{\cos \psi }\sqrt{-1}\left( \frac \partial {\partial \varphi
}+\frac 12\cot \varphi -\frac 12\tan \varphi \right) \gamma _3+\sqrt{-1}%
\left( \frac \partial {\partial \psi }-\frac 32\tan \psi \right) \gamma _4
\end{eqnarray*}
where $\delta _1,\delta _2$ are the derivations on the noncommutative torus
and $\gamma _i$ the Dirac matrices. Using an Ansatz 
\[
f=\sum_{n,m}a_{nm}\left( \varphi ,\psi \right) u^nv^m 
\]
for the Dirac equation 
\[
Df=0 
\]
we get a $q$-deformation of a classical partial differential equation in $%
\varphi $ and $\psi $ for $a_{nm}$. So, solving the Dirac equation on $%
S^{4,\theta }$ translates into a problem of solving a differential equation
for certain $q$-special functions.

From a slightly different perspective, one can also understand this from the
following consideration: It is one of the basic properties of quantum field
theory that one can represent it in terms of infinitely many coupled
harmonic oscillators. This property is not independent of the fact that the
symmetry group of flat three dimensional space is related to $SU\left(
2\right) $ which in turn links to the oscillator algebra. There is a similar
situation for the case of the $S^{4,\theta }$. The algebras $S^{4,\theta }$
contain the suspension of a noncommutative 3-sphere given by an analytic
continuation of the quantum group $SU_q\left( 2\right) $ to the case where $%
q $ is a root of unity. Now, the quantum group $SU_q\left( 2\right) $
relates through the Jordan-Schwinger isomorphism between the (roughly
speaking dual)\ quantum algebra $\mbox{\v{U}}_q(sl_2)$ and the extended
algebra of two $q$-oscillators $\mathcal{A}_q^{ext,2}$ to \bigskip $q$%
-oscillators (see \cite{KS}). By the $SU_q\left( 2\right) $ symmetry of the
suspended noncommutative 3-sphere, it follows that a projective module over $%
S^{4,\theta }$ carries always a coaction of basically $SU_q\left( 2\right) $
on it (taking the coaction of $SU_q\left( 2\right) $ on the suspended
noncommutative 3-sphere and the trivial coaction on the generator $t$ of $%
S^{4,\theta }$). Remember that a projective module is the noncommutative
analog of a vector bundle. Therefore, the quantum state space of a
noncommutative gauge theory over $S^{4,\theta }$ should always carry a
corepresentation of $SU_q\left( 2\right) $. But then, by the results cited
above, such a noncommutative quantum field theory can always be interpreted
as a system of coupled $q$-oscillators. Now, $q$-oscillators relate to $q$%
-Hermite polynomials as $q$-special functions (see \cite{KS}). So, at least
if we treat the quantum field theories on $S^{4,\theta }$ perturbatively,
starting from a system of free $q$-oscillators, we should, again, be led to
expressions in terms of $q$-special functions.

In conclusion, both arguments suggest that perturbative quantum field theory
on $S^{4,\theta }$ should have a formulation in terms of $q$-special
functions. Since in usual quantum gauge theories on commutative space-time,
the weights of Feynman graphs can be calculated as values of so called
multiple polylogarithms (see \cite{BK}, \cite{Kre}) and this fact plays a
considerable role for the structural properties of quantum field theory (see 
\cite{CK}, \cite{Kon}, \cite{KoSo}), a first step toward perturbative
quantum field theory on $S^{4,\theta }$ should be to look for $q$%
-deformations of multiple polylogarithms. For the logarithm and dilogarithm
cases, $q$-deformations have been studied (see e.g. \cite{Kir}, \cite{Koo}
and the literature cited therein) but even for the dilogarithm case, we use
a slightly different form for the $q$-deformation, here. The main advantage
of the form which we give, here, is that the $q$-deformed multiple
polylogarithms form an algebra, then, and for the restriction to the case of 
$q$-deformed multiple $\zeta $-values we even get a noncommutative
deformation of the classical Hopf algebra structure of multiple
polylogarithms (see \cite{Gon}). Since precisely these structural features
of multiple polylogarithms give the above mentioned link to quantum field
theory, we consider this to be an important property for the purpose we have
in mind.

For the reader who is not aquainted with the notations and conventions of $q$%
-calculus and $q$-special functions, we refer e.g. to \cite{KS} for a short
introduction.

\bigskip

\section{$q$-deformed multiple polylogarithms}

Let 
\[
Li_{n_1,...,n_m}\left( z_1,...,z_m;q\right) =\sum_{0<k_1<...<k_m}\frac{%
z_1^{k_1}...z_m^{k_m}}{\left( 1-q^{k_1}\right) ^{n_1}...\left(
1-q^{k_m}\right) ^{n_m}} 
\]
be the $q$-deformed multiple polylogarithms. We have 
\[
\left( 1-q\right) ^{n_1+...+n_m}Li_{n_1,...,n_m}\left( z_1,...,z_m;q\right) 
\stackrel{q\rightarrow 1}{\rightarrow }Li_{n_1,...,n_m}\left(
z_1,...,z_m\right) 
\]
where 
\[
Li_{n_1,...,n_m}\left( z_1,...,z_m\right) =\sum_{0<k_1<...<k_m}\frac{%
z_1^{k_1}...z_m^{k_m}}{k_1^{n_1}...k_m^{n_m}} 
\]
are the classical multiple polylogarithms. With the $q$-derivative defined
as (see e.g. \cite{KS}) 
\[
D_{z,q}f\left( z_0\right) =\frac{f\left( qz_o\right) -f\left( z_o\right) }{%
\left( q-1\right) z_0} 
\]
we have the following lemma:

\bigskip

\begin{lemma}
We have 
\[
D_{z_j,q}Li_{n_1,...,n_m}\left( z_1,...,z_m;q\right) =\frac 1{\left(
1-q\right) z_j}Li_{n_1,...,n_j-1,...n_m}\left( z_1,...,z_m;q\right) 
\]
Especially 
\[
D_{z,q}Li_n\left( z,q\right) =\frac 1{\left( 1-q\right) z}\cdot \frac 1{1-z}
\]
\proof%
By calculation.%
\endproof%
\end{lemma}

\bigskip

In contrast to the $q$-deformations introduced for the logarithm and the
dilogarithm case in \cite{Kir} and \cite{Koo}, we have chosen a version for
the $q$-deformation which has the following important property:

\bigskip

\begin{lemma}
The rational linear span of the $Li_{n_1,...,n_m}\left( z_1,...,z_m;q\right) 
$ forms an associative algebra over $\Bbb{Q}$.
\end{lemma}

\proof%
Follows from the definition.%
\endproof%

\bigskip

A special case of this is the fact that the so called symmetry relations
(see \cite{Gon 1997}) hold also in the $q$-deformed case, as is seen from
the following corollary:

\bigskip

\begin{corollary}
We have 
\[
Li_1\left( x,q\right) Li_1\left( y,q\right) =Li_{1,1}\left( x,y,q\right)
+Li\left( y,x,q\right) +Li_2\left( xy,q\right) 
\]
\end{corollary}

\bigskip

Also there is a $q$-deformed version of the so called distribution relations
(see \cite{Gon 1997} for the classical case. again), now:

\bigskip

\begin{lemma}
We have 
\[
Li_2\left( x,q^n\right) =\frac 1n\sum_{y^n=x}Li_2\left( y,q\right) 
\]
and 
\[
Li_{1,1}\left( x_1,x_2,q^n\right) =\frac 1n\sum_{y_i^n=x_i}Li_{1,1}\left(
y_1,y_2,q\right) 
\]
\end{lemma}

\proof%
Direct calculation.%
\endproof%

\bigskip

An important consequence of the property 
\[
\frac d{dz_j}Li_{n_1,...,n_m}\left( z_1,...,z_m\right) =\frac
1{z_j}Li_{n_1,...,n_j-1,...n_m}\left( z_1,...,z_m\right) 
\]
of the classical multiple polylogarithms is the integral formula 
\[
Li_{n_1,...,n_m}\left( z_1,...,z_m\right) =\int_0^{z_j}\frac
1{w_j}Li_{n_1,...,n_j-1,...,n_m}\left( z_1,...w_j,...,z_m\right) 
\]
which leads as a special case to the integral formula 
\begin{eqnarray*}
\zeta \left( n_1,...,n_m\right) &=&\sum_{0<k_1<...k_m}\frac
1{k_1^{n_1}...k_m^{n_m}} \\
&=&\int_0^1\stackunder{n_1times}{\underbrace{\frac{dt}{1-t}\circ \frac{dt}%
t\circ ...\circ \frac{dt}t}}\circ ...\circ \stackunder{n_mtimes}{\underbrace{%
\frac{dt}{1-t}\circ \frac{dt}t\circ ...\circ \frac{dt}t}}
\end{eqnarray*}
noted by Kontsevich (see e.g. \cite{Gon} for the definition of the right
hand integral) for the multiple $\zeta $-values (i.e. the multiple
polylogarithms with all variables identical to unity and $n_m>1$). E.g. 
\[
\zeta \left( 2\right) =\int_{0\leq t_1\leq t_2\leq 1}\frac{dt_1}{1-t_1}%
\wedge \frac{dt_2}{t_2} 
\]
As a consequence of Lemma 1, we have a similar integral formula for the case
of $q$-deformed multiple $\zeta $-values with the integral replaced by the $%
q $-integral (see e.g. \cite{KS}) and the appropriate powers of $\left(
1-q\right) $ inserted. E.g. 
\[
\zeta _q\left( 2\right) =\frac 1{1-q}\int_{0\leq t_1\leq t_2\leq 1}\frac{%
d_qt_1}{1-t_1}\wedge \frac{d_qt_2}{t_2} 
\]
Up to now, we have considered the case where the multiple polylogarithms
become $q$-deformed but the variables $z_j$ remain classical commutative
numbers from $\Bbb{C}$, only. In the next section, we are going to consider
variables which are $q$-commutative.

\bigskip

\section{$q$-deformed multiple polylogarithms with $q$-commuting variables}

Let $z_j,\ j\in \Bbb{N}$ be $q$-commuting variables, i.e. 
\[
z_iz_j=q\ z_jz_i,\ i<j 
\]
We now have to give an ordering for the variables $z_j$ in the definition of
the $q$-deformed multiple polylogarithms as given above, in order to make
this well defined in the case of $q$-commuting variables. We do this by
requiring that we order the variables from left to right according to
increasing index, e.g. 
\[
Li_{1,1}\left( z_1,z_2;q\right) =\sum_{0<j<k}\frac{z_1^jz_2^k}{\left(
1-q^j\right) \left( 1-q^k\right) } 
\]
but 
\[
Li_{1,1}\left( z_2,z_1;q\right) =\sum_{0<j<k}\frac{z_1^kz_2^j}{\left(
1-q^j\right) \left( 1-q^k\right) } 
\]
With this definition, we have, again, the following property:

\bigskip

\begin{lemma}
For the ordered definition of the $q$-deformed multiple polylogarithms with $%
q$-commuting variables as given above, the rational span of the $q$-deformed
multiple polylogarithms forms an algebra over $\Bbb{Q}$. Especially, the
symmetry relations hold also for the $q$-commuting variables, then.
\end{lemma}

\proof%
Direct calculation, again.%
\endproof%

\bigskip

Especially, this means that the commutator of two $q$-deformed multiple
polylogarithms can be expressed as a sum of such.

Since the differential equation for $q$-deformed multiple polylogarithms
remains formally valid in the case of $q$-commuting variables, too, we can
next ask for an integral formula for the $q$-deformed multiple $\zeta $%
-values in this case. In trying to derive such a formula, one makes a
strange observation: On the one hand, we have 
\[
\zeta _q\left( n_1,...,n_m\right) =\sum_{0<k_1<...k_m}\frac 1{\left(
1-q^{k_1}\right) ^{n_1}...\left( 1-q^{k_m}\right) ^{n_m}} 
\]
i.e. the $\zeta _q\left( n_1,...,n_m\right) $ do not depend on the variables 
$z_j$ at all and therefore remain commuting with each other. On the other
hand, consider e.g. the case of $\zeta _q\left( 2\right) $ where the
integral formula reads 
\[
\zeta _q\left( 2\right) =\frac 1{1-q}\int_{0\leq t_1\leq t_2\leq 1}\frac{%
d_qt_1}{1-t_1}\wedge \frac{d_qt_2}{t_2} 
\]
Now, in a formal sense one would, in the case of $q$-commuting variables,
consider the integrals not to satisfy the Fubini theorem if one multiplies
to expressions of this kind but one would expect the differentials formally
to exchange with a factor $q^{-1}$. So, formally this would lead to an
exchange rule 
\[
\zeta \left( n_{1,}...,n_m;q\right) \zeta \left( j_1,...,j_l;q\right)
=q^{\left( \sum n_i\right) \left( \sum j_i\right) }\zeta \left(
j_1,...,j_l;q\right) \zeta \left( n_{1,}...,n_m;q\right) 
\]
for the $q$-deformed multiple $\zeta $-values for $n_1<...<n_m<j_1<...<j_l$
in the ordering for the $q$-commutativity relations. This means that we
should allow the arguments of the $q$-deformed multiple $\zeta $-values -
i.e. the indices of the $q$-deformed multiple polylogarithms - to become
noncommuting variabels, too. In the sequel, we will assume this to be the
case and the above relation for the $q$-deformed multiple $\zeta $-values to
hold true.

\bigskip

We next come to the question of Hopf algebra structure. On the algebra of
classical multiple polylogarithms there is a coproduct given (see \cite{Gon}
for the definition). Since the coproduct $\Delta $ preserves the subalgebra
of classical multiple $\zeta $-values, we get a bialgebra $\mathcal{Z}$ of
classical multiple $\zeta $-values in this way. For the algebra $\mathcal{Z}%
_q$ of $q$-deformed multiple $\zeta $-values as just defined, $\Delta $ is
no longer compatible as a coproduct. The reason for this is that the number
of arguments - which determines the power of $q$ in the above $q$%
-commutation relation - varies between the different terms of the coproduct
(see the defintion in \cite{Gon}). But this observation immediately suggests
also the solution: One has to formally insert an appropriate number of $q$%
-commuting variables into the different terms of the coproduct in order to
make it an algebra morphism of the deformed product. In conclusion, we have
observed the following result, therefore:

\bigskip

\begin{lemma}
There is a $q$-deformation $\Delta _q$ of the coproduct $\Delta $ of $%
\mathcal{Z}$ which turns $\mathcal{Z}_q$ into a bialgebra. Since there is a
full Hopf algebra structure for $\mathcal{Z}$ (see \cite{Gon}), we can by a
general result in Hopf algebra cohomology (see e.g \cite{CP}) assume without
loss of generality that $\mathcal{Z}_q$ carries a full Hopf algebra
structure, too.
\end{lemma}

\bigskip

\begin{remark}
The commutators for the full algebra of $q$-deformed multiple polylogarithms
- with $q$-commutation relations for both, the variables and the indices -
are much too complicated to allow at present for an answer to the question
if there exists a deformation of the full Hopf algebra of classical multiple
polylogarithms. But since the Hopf structure is intimately tied to the
motivic origin of these periods (see \cite{Gon}), we expect that as a
consequence of ultrarigidity (see \cite{Sch 2001a}, \cite{Sch 2001b}), there
exists no such deformation (because the multiple $\zeta $-values are
supposed to give precisely those periods which are linked to the
Grothendieck-Teichm\"{u}ller group, see \cite{Kon}).
\end{remark}

\begin{remark}
The Hopf algebra of Connes and Kreimer is supposed to be related to the Hopf
algebra of classical multiple $\zeta $-values (see \cite{CK}, \cite{CK 1999}%
). So, we expect $\mathcal{Z}_q$ to be related to a $q$-deformation of the
Connes-Kreimer algebra (see also \cite{GS} for the question of deformations
of the Connes-Kreimer algebra).
\end{remark}

\bigskip

\begin{conjecture}
By the arguments given in the introductory section and by the connection
between the Hopf algebra structure of multiple $\zeta $-values and the
Connes-Kreimer algebra mentioned in the previous remark, one expects $%
\mathcal{Z}_q$ to act on quantum field theories on $S^{4,\theta }$. On the
other hand, we know that the Hopf algebra $\widetilde{\mathcal{H}_{GT}}$,
introduced in \cite{Sch 2001a} as a quantum analog of the
Grothendieck-Teichm\"{u}ller group, acts on quantum field theories on $%
S^{4,\theta }$ (see \cite{Sch 2001b}). We conjecture that the corresponding
representations of $\mathcal{Z}_q$ and $\widetilde{\mathcal{H}_{GT}}$,
gained in this way, are isomorphic.
\end{conjecture}

\bigskip

In conclusion, we have seen that it is possible to introduce a special $q$%
-deformation of multiple polylogarithms which allows for quantum
counterparts of several of the algebraic features of periods which have
become increasingly important in the study of quantum field theories in
recent years.

\bigskip

\textbf{Acknowledgements:}

I thank H. Grosse for discussions on the topics involved and A. Goncharov
for very helpful explanations on classical motivic structures and for
drawing my attention to an important reference. Besides this, I thank the
Deutsche Forschungsgemeinschaft (DFG) for support by a research grant and
the Erwin Schr\"{o}dinger Institute for Mathematical Physics, Vienna, for
hospitality.

\end{document}